\documentclass{IEEEtran}
\usepackage{cite}
\usepackage{amsmath,amssymb,amsfonts}
\usepackage{algorithm}
\usepackage{graphicx}
\usepackage{textcomp}
\usepackage{algpseudocode}
\usepackage{mathrsfs}
\usepackage{mathtools}
\usepackage{scalerel}
\usepackage{setspace}
\usepackage{tikz}
\usepackage{pgfplots}
\usepackage[hidelinks]{hyperref}
\usepgfplotslibrary{groupplots}
\pgfplotsset{compat=1.18}

\mathtoolsset{showonlyrefs}

\everypar{\looseness=-1}
\allowdisplaybreaks
\IEEEaftertitletext{\vspace{-2\baselineskip}}

\def\BibTeX{{\rm B\kern-.05em{\sc i\kern-.025em b}\kern-.08em
    T\kern-.1667em\lower.7ex\hbox{E}\kern-.125emX}}
    
\begin{document}
\title{A Parallel-in-Time Newton's Method for \\ Nonlinear Model Predictive Control}
\author{Casian Iacob, Hany Abdulsamad, Simo S\"arkk\"a%
\thanks{This work was supported by the Business Finland project 9585/31/2021 (AIMODE) and the Finnish Center for Artificial Intelligence (FCAI). (\textit{Correspondng author: Casian Iacob})}
\thanks{The authors are with the Electrical Engineering and Automation department of Aalto University, Espoo, FI-00076 Finland (e-mail: casian.iacob@aalto.fi, hany.abdulsamad@aalto.fi, simo.sarkka@aalto.fi). }
}

\maketitle

\begin{abstract}
 Model predictive control (MPC) is a powerful framework for optimal control of dynamical systems. However, MPC solvers suffer from a high computational burden that restricts their application to systems with low sampling frequencies. This issue is further amplified in nonlinear and constrained systems that require nesting MPC solvers within iterative procedures. In this paper, we address these issues by developing parallel-in-time algorithms for constrained nonlinear optimization problems that take advantage of massively parallel hardware to achieve logarithmic computational time scaling over the planning horizon. We develop time-parallel second-order solvers based on interior point methods and the alternating direction method of multipliers, leveraging fast convergence and lower computational cost per iteration. The parallelization is based on a reformulation of the subproblems in terms of associative operations that can be parallelized using the associative scan algorithm. We validate our approach on numerical examples of nonlinear and constrained dynamical systems. 
\end{abstract}

\setlength{\textfloatsep}{0pt}
\setlength{\intextsep}{0pt}
\section{Introduction}
Nonlinear optimal control problems \cite{lewis2012optimal} are inherent in advanced engineering systems in robotics and automotive applications. Model predictive control (MPC) \cite{rawlings2017model} provides a principled framework for addressing these problems through numerical optimization, particularly in constrained settings. In this paper, the aim is to develop parallel-in-time methods for nonlinear constrained MPC problems.

Optimal control algorithms can be classified based on either their problem formulation or the optimization techniques they employ. The overall problem can be stated and solved either as a condensed nonlinear program or recursively via dynamic programming \cite{borrelli2017predictive, rawlings2017model}. A further classification arises from the choice of the optimization method. In this paper, we focus on constrained nonlinear optimization problems where suitable solvers are, for example, the alternating direction method of multipliers (ADMM) \cite{boyd2011distributed} and interior point (IP) methods \cite{nocedal1999numerical}.

The different optimization methods employed in condensed and dynamic programming approaches for MPC offer distinct advantages in computational efficiency. In linear model predictive control, condensation techniques construct a linearly-constrained quadratic program that can be solved either by ADMM \cite{boyd2011distributed, stellato2020osqp} or IP \cite{rao1998application}. For nonlinear MPC, condensation IP methods have been extended to general nonlinear settings based on sequential quadratic programming \cite{rawlings2017model}. Recently, dynamic programming adaptations of IP \cite{chen2019autonomous, pavlov2021interior} and ADMM \cite{ma2022alternating} have gained popularity in robotics. These methods leverage the sparse structure of the problem by using differential dynamic programming (DDP) \cite{jacobson1970differential} or the iterative linear quadratic regulator (iLQR) \cite{tassa2012synthesis}. The advantage of these algorithms lies in their explicit recursive form, which results in linear complexity with respect to the planning horizon, independent of the software framework used. In contrast, condensation solutions typically require leveraging explicit matrix factorization techniques to achieve comparable computational efficiency.

Despite its widespread interest and practical applications, Model Predictive Control (MPC) faces significant challenges in real-world deployment. One of the primary obstacles is the high computational complexity associated with solving control problems efficiently \cite{rawlings2017model}. The computational burden can be overwhelming since MPC requires continuous re-planning within short time intervals dictated by system dynamics. While recursive algorithms help reduce complexity to linear time, they often fall short of delivering the necessary speed for real-time applications. Consequently, MPC is most commonly used in systems with slower response times, where computational demands are more manageable.

Recent advancements in parallel computing offer potential solutions to overcome the computational challenges of MPC, especially for high-frequency systems. Optimization methods tailored for parallel execution, such as those in \cite{schubiger2020gpu}, \cite{bishop2024relu}, and, \cite{persson2024optimization}, focus on constructing convex quadratic programs commonly used in condensed linear MPC, while \cite{pacaud2024gpu}, focuses on solving the nonlinear Karush--Kuhn--Tucker system arising in the primal-dual approach. These methods leverage hardware-optimized factorization algorithms for efficient matrix-vector operations. A key limitation of these methods is their dependence on the performance of sparse matrix parallelization routines, lacking a more explicit computational structure. 

The contribution of this paper is to develop computationally efficient nonlinear MPC solvers, by temporal parallelization of prefix-sum operations, inspired by the methods in \cite{sarkka2022temporal}. Temporal parallelization speeds up dynamic programming by eliminating recursive dependencies, reducing computational time, and achieving logarithmic span complexity relative to the planning horizon. The appeal of these approaches lies in their explicit algorithmic design, which can be implemented on a wide range of software and hardware architectures. We integrate temporal parallelization in an explicit dynamic programming implementation of Newton's method for optimal control \cite{dunn1989efficient, de1988differential, roulet2022iterative}, and combine it with an IP method and ADMM. We implement both solvers on GPUs and demonstrate their computational efficiency and scalability on a set of numerical examples, validating the computational advantage of our work. The faster runtimes are obtained for problems with long optimization horizons. The computational improvement is negligible for optimization over short horizons.  

The remainder of the paper is structured as follows. In Section~II, we present the constrained nonlinear optimization problem underpinning nonlinear model predictive control and outline the primal log-barrier method and the alternating direction method of multipliers, followed by an introduction to prefix-sum operations and associative scans. In Section~III, we present our parallel-in-time Newton solver that can be adapted for both optimization techniques. Finally, in Section~IV, we evaluate our approach.

\section{Background}
\label{sec:background}
In this section, we start with an overview of IP methods \cite{nocedal1999numerical, chen2019autonomous} and ADMM \cite{boyd2011distributed, ma2022alternating} in the context of constrained nonlinear optimal control. Then, we give a brief introduction to parallel computation for prefix-sums operations \cite{blelloch1990prefix}.

\subsection{Constrained Nonlinear Optimal Control Problem} \label{sec:constrained-optimal-control}
In this paper, we consider the following constrained optimal control problem:
\begin{equation}
\begin{split}
    \label{eq:ocp}
     &\min_{u_{1:N}} \,\,\underbrace{l_{N+1}(x_{N+1}) + \sum_{t=1}^{N} l_{t}(x_{t}, u_{t})}_{\coloneqq J(u_{1:N})} 
     \\ &\textrm{subject\,to} \\
    &\begin{aligned}
        & x_{t+1} = f_{t}(x_{t}, u_{t}), && \forall \, t=1, \hdots, N, \\
        & g_{t}(x_{t}) \leq 0, && \forall \, t=1, \hdots, N, \\
        & h_{t}(u_{t}) \leq 0, && \forall \, t=1, \hdots, N, \\
        & x_{1} = \hat{x}_{1}, &&
    \end{aligned}
\end{split}
\end{equation}
where $x_t \in \mathbb{R}^{d_x}$ is the state vector, $u_t \in \mathbb{R}^{d_u}$ is the control vector, and $f_{t}:\mathbb{R}^{d_x} \times \mathbb{R}^{d_u} \rightarrow\mathbb{R}^{d_x}$ is the dynamic model function. The stage and terminal costs are $l_{t}: \mathbb{R}^{d_x} \times \mathbb{R}^{d_u} \rightarrow \mathbb{R}$ and $l_{N+1}: \mathbb{R}^{d_x} \rightarrow \mathbb{R}$, respectively. 
The functions $g_{t}: \mathbb{R}^{d_x} \rightarrow \mathbb{R}^{d_x}$ and $h_{t}: \mathbb{R}^{d_u} \rightarrow \mathbb{R}^{d_u}$ enforce inequality constraints on state and control.

Next, we review two types of methods for solving these problems: primal log-barrier IP methods and ADMM.

\subsection{Primal Log-Barrier Interior Point Method}
The primal log-barrier approach is part of the IP methods family \cite{nocedal1999numerical} which minimizes a cost function augmented with the log-barrier function of the constraints. We define the barrier problem corresponding to \eqref{eq:ocp} as
\begin{equation}
    \begin{split}
        &\min_{u_{1:N}}\,\, \underbrace{
        \begin{aligned}[t]
            J(u_{1:N}) & - \mu \, \sum_{t=1}^{N} \sum_{m=1}^{d_x} \log(-g_{t}^{m}(x_{t})) \\
            \qquad & - \mu \, \sum_{t=1}^{N} \sum_{m=1}^{d_u} \log(-h_{t}^{m}(u_{t}))
        \end{aligned}
        }_{\coloneqq J_{\mathrm{LB}}(u_{1:N},\mu)}\\
    &\textrm{subject\,\,to}\\
    &x_{t+1} = f_t(x_t, u_t), \quad x_{1} = \hat{x}_{1},
    \end{split}
\end{equation}
where $\mu > 0$ is the barrier parameter. By sequentially solving barrier problems and decreasing the barrier parameter, the constrained solution to \eqref{eq:ocp}  is recovered \cite{nocedal1999numerical}. The idea is that as $\mu \rightarrow 0$, the log-barrier function becomes an indicator function. The indicator function $\mathcal{I}_\mathcal{C}(x)$ for a general set $\mathcal{C}$ is defined as 
\begin{equation}
    \mathcal{I}_\mathcal{C}(x) = \begin{cases}0, & x\in\mathcal{C}\\ \infty, & x\notin{\mathcal{C}}\end{cases}.
\end{equation}
The primal log-barrier method is summarized as:
\begin{align}
    u_{1:N}^{(k+1)} & \gets \arg\min_{{u_{1:N}}} J_\mathrm{LB}(u_{1:N}, \mu^{(k)}), \label{eq:lb-min} \\
    \mu^{(k+1)} & \gets \zeta \, \mu^{(k)}, \label{eq:lb-step-2}
\end{align}
where $k$ is the iteration index and $0 < \zeta < 1$ is a step size. Steps \eqref{eq:lb-min} and \eqref{eq:lb-step-2} are repeated until the barrier parameter $\mu$ converges to a predefined tolerance. 

\subsection{Alternating Direction Method of Multipliers}
The alternating direction method of multipliers \cite{boyd2011distributed} is a distributed version of the augmented Lagrangian method of multipliers where the objective is split between the original decision variable and a constrained consensus variable. Let us denote the concatenated vector of inequality constraints by $w(x_t, u_t) = \begin{bmatrix}g_t^\top(x_t) & h_t^\top(u_t)\end{bmatrix}^\top$ and the sets defined by the dynamic and inequality constraints as
\begin{equation}
    \begin{split}
        \mathcal{F} &= \{(x_t, u_t) \mid x_{t+1} = f_t(x_t, u_t),\,x_1 = \hat{x}_1\},\\
        \mathcal{C} &= \{(x_t, u_t) \mid w(x_t, u_t) \leq 0\}.
    \end{split}
\end{equation}

In ADMM, to enforce the constraints, we introduce a consensus variable $z_t \in \mathbb{R}^{d_x + d_u}$ such that $z_t = w_{t}(x_{t}, u_{t})$. The problem defined in \eqref{eq:ocp} can now be rewritten as
\begin{equation}
    \begin{split}
        &\min_{u_{1:N}, z_{1:N}} J(u_{1:N}) + \mathcal{I}_{\mathcal{F}}(x_{1:N+1}, u_{1:N}) + \mathcal{I}_{\mathcal{C}}(z_{1:N})\\
        &\textrm{subject to}\\
        &z_{t} = w_{t}(x_{t}, u_{t}), \quad \forall \, t=1, \hdots, N,
    \end{split}
\end{equation}
where $\mathcal{I}_{\mathcal{F}}$ and $\mathcal{I}_\mathcal{C}$ are indicator functions for the sets $\mathcal{F}$ and $\mathcal{C}$, respectively. The above optimization problem is tackled by constructing the augmented Lagrangian as 
\begin{equation}
    \begin{split}
        &\mathcal{L}_{\rho}(x_{1:N+1}, u_{1:N}, z_{1:N}, v_{1:N})\\
        &=
        \underbrace{J(u_{1:N}) + \frac{\rho}{2} \sum_{t=1}^{N} \left\Vert w_{t}(x_t, u_{t}) - z_{t} + \rho^{-1} v_{t} \right\Vert_{2}^{2}}_{\coloneqq J_{\mathrm{AD}}(u_{1:N}, z_{1:N}, v_{1:N})} \\
        &\qquad + \frac{1}{2 \rho} \sum_{t=1}^{N} \left\Vert v_{t} \right\Vert_{2}^{2} + \mathcal{I}_{\mathcal{F}}(x_{1:N}, u_{1:N}) + \mathcal{I}_{\mathcal{C}}(z_{1:N}),
    \end{split}
\end{equation}
where $\rho > 0$ is a penalty parameter, and $v_t \in \mathbb{R}^{d_x+d_u}$ are Lagrange multipliers. Optimizing this Lagrangian with respect to $u_{1:N}$, $z_{1:N}$, and $v_{1:N}$ leads to the ADMM steps:
\begin{align}
    \begin{split}
        u_{1:N}^{(k+1)} &\gets
         \arg\min_{{u_{1:N}}} J_{\mathrm{AD}}(u_{1:N}, z_{1:N}^{(k)}, v_{1:N}^{(k)})\\
         &\textrm{subject to}\\
         &x_{t+1} = f_t(x_t, u_t), \,\, x_{1} = \hat{x}_{1}
    \end{split}\label{eq:admm-min}\\
    x_{t+1}^{(k+1)} & \gets f_t(x_t^{(k+1)}, u_{t}^{(k+1)}), \label{eq:admm-step-2}\\
    z_{t}^{(k+1)} & \gets \Pi_{\mathcal{C}}(w_t(x_{t}^{(k+1)}, u_{t}^{(k+1)}) + \rho^{-1} v_{t}^{(k)}), \label{eq:admm-step-4}\\
    v_{t}^{(k+1)} & \gets v_{t}^{(k)} + \rho \, (y_{t}^{(k+1)} - z_{t}^{(k+1)}), \label{eq:admm-step-5}
\end{align}
$\forall \, t=1, \hdots, N$, where $k$ is the iteration index and $\Pi_{\mathcal{C}}$ is the Euclidean projection operator on the inequality constraints. The efficiency of ADMM depends heavily on the solution to the projection step \eqref{eq:admm-step-4} that should be available in closed form. For some special cases of projections of a point on a convex set, see \cite{boyd2004convex}. The steps \eqref{eq:admm-min}, \eqref{eq:admm-step-2}, \eqref{eq:admm-step-4}, and \eqref{eq:admm-step-5} are repeated until the follow primal and dual residuals converge:
\begin{equation*}
   r_p^{(k)} = w_{1:N}(x_{1:N}^{(k)}, u_{1:N}^{(k)}) - z_{1:N}^{(k)}, \quad r_d^{(k)} = z_{1:N}^{(k)} - z_{1:N}^{(k-1)},
\end{equation*}

We observe that IP and ADMM methods require repeatedly solving equality-constrained nonlinear optimal control sub-problems defined in \eqref{eq:lb-min} and \eqref{eq:admm-min}, respectively. Furthermore, due to the non-convexity of these sub-problems, each instance requires multiple internal iterations to converge. If the number of iterations of IP or ADMM is $K$ and each sub-problem requires $M$ internal iterations, this leads to a total of $K M$ iterations. Therefore, it is crucial to develop solvers that converge fast, minimizing the number of iterations, and are also computationally efficient per sub-problem. To address this issue, we propose a massively-parallelizable Newton-type solver for \eqref{eq:lb-min} and \eqref{eq:admm-min}, which enjoys favorable convergence and computational properties. 

\subsection{Parallel Computation via Prefix-Sum Operations}
\label{sec:prefix-sum}
The constrained nonlinear optimal control sub-problems from Section~\ref{sec:constrained-optimal-control} can be tackled via methods such as DDP \cite{jacobson1970differential} or iLQR \cite{tassa2012synthesis}. These methods perform recursive dynamic programming computations with linear time complexity $\mathcal{O}(N)$ per iteration with respect to the planning horizon $N$.

In \cite{sarkka2022temporal}, the authors propose a parallel-in-time realization of general dynamic programming equations, including Riccati recursions based on a prefix-sum computation. Such computations can be formulated as a so-called parallel \emph{prefix-sum} or \emph{scan} operation \cite{blelloch1990prefix}. Given a sequence of elements $[a_1, a_1, \hdots, a_{N}]$ and an associative operator $\otimes$ the prefix-sum operation returns a new sequence of the form $[a_1, (a_1 \otimes a_2), \hdots, (a_1 \otimes a_2 \hdots \otimes a_{N})]$. By a suitable definition of the elements and the associative operator an optimal control problem can be solved via the parallel scan algorithm \cite{sarkka2022temporal} which has $\mathcal{O}(\log N)$ time complexity, drastically improving the computational efficiency of recursive computation in dynamic programming algorithms. 

However, the Riccati-equation-based approach in \cite{sarkka2022temporal} is limited to linear or \emph{linearized} dynamics, corresponding to an iLQR-style iterated Gauss-Newton method. For constrained, non-convex problems such as those discussed in Section~\ref{sec:constrained-optimal-control}, it may be advantageous to leverage second-order information and perform full Newton steps to improve convergence \cite{nocedal1999numerical} and reduce the overall number of iterations. Moreover, while DDP-style methods enjoy quadratic convergence, they do not appear to be amenable to the same kind of parallelization. This is due to additional terms that break the associativity of the proposed operator, see Lemma~10 in \cite{sarkka2022temporal}.

To address this gap, we present an alternative parallel-in-time Newton-type solver for nonlinear optimal control problems with an iterative procedure over three associative scans. We provide the definitions of corresponding associative operators $\otimes$ and elements $a_{1:N}$. We demonstrate that this parallel solver leads to a significant computational advantage over equivalent sequential solutions.

\section{Parallel Dynamic Programming Realization of Newton's Method}
\label{sec:parallel-newton}
In this section, we rewrite sub-problems \eqref{eq:lb-min} and \eqref{eq:admm-min} in a general form. We then introduce a parallel-in-time Newton's method for optimal control by adapting the steps of \cite{dunn1989efficient}. The method is divided into two backward passes and one forward pass. We define associative elements and an associative operation for each of these sequential operations such that the computation is performed in parallel via associative scans.

\subsection{Augmented Costs and Lagrangians for IP and ADMM}
The minimization steps \eqref{eq:lb-min} and \eqref{eq:admm-min} are instances of nonlinear optimal control problems \cite{lewis2012optimal}. We can write these objectives as the original cumulative cost augmented by a barrier or penalty function, respectively. Therefore we consider the following augmented objective:
\begin{equation}
    J_{c}(u_{1:N}) = l_{N+1}(x_{N+1}) + \sum_{t=1}^{N} \Big[ l_t(x_t, u_t) + c_t(x_t, u_t) \Big],
\end{equation}
where the augmenting terms for interior point methods are
\begin{equation}
    \hat{c}_t(x_t, u_t) = - \mu \left[ \sum_{m=1}^{d_x} \log(-g_t^m(x_{t})) + \sum_{m=1}^{d_u} \log(-h_t^m(u_{t})) \right].
\end{equation}
In the case of ADMM, we have
\begin{equation}
    \bar{c}_t(x_t, u_t) = \frac{\rho}{2} \left\Vert w_t(x_{t}, u_{t}) - z_{t} + \rho^{-1} v_t \right\Vert_2^2.
\end{equation}
Thus we can rewrite the problems \eqref{eq:lb-min} and \eqref{eq:admm-min} as
\begin{equation}
    \label{eq:augmented-ocp}
    \begin{split}
        &\min_{u_{1:N}} \,\, l_{N+1}(x_{N+1}) + \sum_{t=1}^{N} \Big[ l_t(x_t, u_t) + c_t(x_t, u_t) \Big],\\
        &\textrm{subject to}\\
        &x_{t+1} = f_t(x_t, u_t), \quad x_{1} = \hat{x}_{1},
    \end{split}
\end{equation}
where $c_t = \hat{c}_t$ or $c_t = \bar{c}_t$. To incorporate the dynamics constraint, we formulate the Lagrangian of the problem \eqref{eq:augmented-ocp}
\begin{equation}
    \label{eq:ocp-lagrangian}
    \begin{split}
        \mathcal{L}_{c}(u_{1:N}, \lambda_{1:N+1}) = {} & l_{N+1}(x_{N+1}) + \lambda_{1}^{\top} \hat{x}_{1} \\
         &  + \sum_{t=1}^{N} H_t(x_t, u_t) - \sum_{t=1}^{N+1} \lambda_t^{\top} x_t ,
    \end{split}
\end{equation}
where $ H_t(x_t, u_t) = l_t(x_t, u_t) + c_t(x_t, u_t) + \lambda_{t+1}^\top f_t(x_t, u_t)$ is the Hamiltonian and $\lambda_t \in \mathbb{R}^{d_x}$ are the co-states.

We can minimize \eqref{eq:augmented-ocp} with Newton's method following an iterative recursive procedure according to \cite{dunn1989efficient}. At every iteration $(n)$, starting from a control sequence $u_{1:N}^{(n)}$ and a corresponding state trajectory $x_{1:N+1}^{(n)}$, the optimal co-states $\lambda_{1:N+1}^{(n)}$ are computed. Conditioned on these co-states, the Lagrangian function \eqref{eq:ocp-lagrangian} is expanded up to the second order in the neighborhood of the nominal trajectory $\{ x_{1:N+1}^{(n)}, u_{1:N}^{(n)} \}$. This leads to a convex optimization problem over the control deviation vector $\delta u_{1:N} = u_{1:N} - u_{1:N}^{(n)}$. The resulting problem resembles the linear quadratic regulator and is solved by invoking Bellman's principle of optimality \cite{lewis2012optimal, kirk2004optimal}.

To avoid a cluttered notation, in the upcoming sections, we will drop the iteration signifier $(n)$, and instead, we refer to \emph{a} nominal trajectory with $\{ \hat{x}_{1:N+1}, \hat{u}_{1:N} \}$. In the following, we demonstrate that such a Newton-type scheme corresponds to three passes, each fully parallelizable with associative scans. In line with the procedure described in Section~\ref{sec:prefix-sum}, we provide a definition of the operator and corresponding elements for each of these associative scans. 
\subsection{Parallel Co-State Pass}
\label{sec:co-states-scan}
Following \cite{lewis2012optimal}, we can derive the optimal solution for the co-states $\lambda_{t}^*$, $\forall \, 1 \leq t \leq N+1$, by setting the derivative of the Lagrangian in \eqref{eq:ocp-lagrangian} with respect to $x_{t}$ to zero, leading to the following optimality condition:
\begin{equation}
    \label{eq:co-state-eq}
    \lambda_t^* = \frac{\partial H_t}{\partial x_t}(\lambda_{t+1}) = \frac{\partial l_t}{\partial x_t} +\frac{\partial c_t}{\partial x_t} + \Big[ \lambda_{t+1}^* \Big]^{\top} \frac{\partial f_t}{\partial x_t},
\end{equation}
with the following boundary condition at $N+1$:
\begin{equation}
    \label{eq:co-state-final}
    \lambda_{N+1}^* = \frac{\partial l_{N+1}}{\partial x_{N+1}}.
\end{equation}
The co-states are computed via a backward recursion with affine dependencies. We can rewrite this operation using function compositions in terms of the final co-state as follows
\begin{equation}
    \lambda_t^* = \left( \frac{\partial H_t}{\partial x_t} \circ \cdots \circ \frac{\partial H_N}{\partial x_N} \right) (\lambda_{N+1}^*).
\end{equation}
We can select the function composition $\circ$ to be the associative operator in the parallel associative scan algorithm:
\begin{equation}
    \label{eq:co-state-operator}
    \frac{\partial H}{\partial x} \Big\vert_{j,i} = \frac{\partial H}{\partial x} \Big\vert_{j,t} \circ \frac{\partial H}{\partial x} \Big\vert_{t,i},
\end{equation}
where $\cdot\vert_{j,i}$ denotes an associative element resulting from the composition of the Hamiltonian derivative $\partial H / \partial x$ between two time steps $j \rightarrow i$ with $i > j$:
\begin{align}
    \label{eq:co-state-elements}
    \frac{\partial H}{\partial x} \Big\vert_{j,i} & = \left( \frac{\partial H_{j}}{\partial x_{j}} \circ \cdots \circ \frac{\partial H_{i}}{\partial x_{i}} \right) (\lambda_{i+1}^*) \\
    & = \frac{\partial l}{\partial x} \Big\vert_{j,i} + \frac{\partial c}{\partial x} \Big\vert_{j,i} + \frac{\partial f}{\partial x} \Big\vert_{j,i}^{\top} \, \lambda_{i+1}^*.
\end{align}
The associative operator defined in~\eqref{eq:co-state-operator} combines the terms of two elements defined in \eqref{eq:co-state-elements} as follows
\begin{equation}
    \label{eq:co-states-comb-rule}
    \begin{split}
        \frac{\partial l}{\partial x} \Big\vert_{j,i}  & = \frac{\partial l}{\partial x} \Big\vert_{j,t} + \frac{\partial f}{\partial x} \Big\vert_{j,t}^\top \, \frac{\partial l}{\partial x} \Big\vert_{t,i}, \\
        \frac{\partial c}{\partial x} \Big\vert_{j,i} & = \frac{\partial c}{\partial x} \Big\vert_{j,t} + \frac{\partial f}{\partial x}  \Big\vert_{j,t}^\top \, \frac{\partial c}{\partial x} \Big\vert_{t,i}, \\
        \frac{\partial f}{\partial x} \Big\vert_{j,i}^{\top} & = \frac{\partial f}{\partial x} \Big\vert_{j,t}^{\top} \,\, \frac{\partial f}{\partial x} \Big\vert_{t,i}^{\top}.
    \end{split}
\end{equation}
The associative scan should be started from the initial elements defined as
\begin{equation}
    a_{t} = \frac{\partial H}{\partial x} \Big\vert_{t,t+1}
    = \frac{\partial l}{\partial x} \Big\vert_{t,t+1}
    + \frac{\partial c}{\partial x} \Big\vert_{t,t+1}
    + \frac{\partial f}{\partial x} \Big\vert_{t,t+1}, 
\end{equation}
for $1 \leq t \leq N$, where
\begin{equation}
    \label{eq:co-state-elem-init}
    \frac{\partial l}{\partial x} \Big\vert_{t,t+1} = \frac{\partial l_t}{\partial x_t}, \,\,
    \frac{\partial c}{\partial x} \Big\vert_{t,t+1} = \frac{\partial c_t}{\partial x_t}, \,\,
    \frac{\partial f}{\partial x} \Big\vert_{t,t+1} = \frac{\partial f_t}{\partial x_t},
\end{equation}
$\forall \, 1 \leq t < N$, while the element $a_N$ is initialized with
\begin{equation}
    \label{eq:co-state-final-elem-init}
    \begin{split}
    \frac{\partial l}{\partial x} \Big\vert_{N,N+1} & = \frac{\partial l_N}{\partial x_N} + \left[\frac{\partial f_N}{\partial x_N}\right]^{\top} \! \lambda_{N+1}^*, \\
    \frac{\partial c}{\partial x} \Big\vert_{N,N+1} & = \frac{\partial c_N}{\partial x_N},\quad
    \frac{\partial f}{\partial x} \Big\vert_{N,N+1} = 0.
    \end{split}
\end{equation}
Note that the evaluation of the derivatives of $l_t(\cdot), c_t(\cdot)$ and $f_t(\cdot)$ is done with respect to a nominal trajectory $\{ \hat{x}_{1:N+1}, \hat{u}_{1:N} \}$, that corresponds to an intermediate solution in the iterative Newton method.
\begin{algorithm}[t]
    \caption{Parallel Co-State Pass}
    \label{alg:co-state-parallel}
    \begin{algorithmic}[1]
        \State \textbf{Input}: Nominal: $\{ \hat{x}_{1:N+1}, \hat{u}_{1:N} \}$.
        \State Evaluate the partial derivative of $l_{N+1}(\cdot)$ at $\hat{x}_{N+1}$ and compute the final co-state $\lambda_{N+1}^{*}$ from \eqref{eq:co-state-final}.
        \State Evaluate all partial derivatives of $l_{1:N}(\cdot)$, $c_{1:N}(\cdot)$, and $f_{1:N}(\cdot)$ at $\{ \hat{x}_{1:N}, \hat{u}_{1:N} \}$ required for \eqref{eq:co-state-final-elem-init} and \eqref{eq:co-state-elem-init}.
        \State Initialize associative element $a_N$ according to \eqref{eq:co-state-final-elem-init}.
        \State Initialize associative elements $a_{1:N-1}$ according to \eqref{eq:co-state-elem-init}.
        \State Define the associative operator $\otimes$ as in \eqref{eq:co-state-operator}.
        \State Execute a parallel associative scan to compute  $\lambda_{t}^* = a_{t} \otimes \cdots \otimes a_N$, for $t= 1, \ldots, N$.
        \State \textbf{Output}: $\lambda_{1:N+1}^*$.
    \end{algorithmic}
\end{algorithm}
Given these definitions of the associative operator~\eqref{eq:co-state-operator} and its individual elements~\eqref{eq:co-state-elements}, the computation in~\eqref{eq:co-state-eq} can be performed logarithmic time $\mathcal{O}(\log N)$ with the parallel associative scan scheme discussed in Section~\ref{sec:prefix-sum}
\begin{equation}
    \label{eq:co-states-prefix-sum}
     \lambda_{t}^* = a_{t} \otimes \cdots \otimes a_N, \quad \forall \, 1 \leq t \leq N.
\end{equation}
Algorithm~\ref{alg:co-state-parallel} summarizes the parallel co-state pass. 
\subsection{Parallel Value Function Pass}
In this section, we follow the general scheme to derive a recursive Newton step as described in \cite{dunn1989efficient} and show that the resulting dynamic programming solution can be described in terms of a parallelizable associative scan. Given a nominal control sequence $\hat{u}_{1:N}$ and the respective co-states $\lambda_{1:N+1}^{*}$, we continue by expanding the Lagrangian~\eqref{eq:ocp-lagrangian} with respect to $u_{1:N}$ up to the second order 
\begin{equation}
    \label{eq:quad-lagrangian}
    \mathcal{L}_{c}(u_{1:N}) \approx 
    \begin{aligned}[t]
        & \mathcal{L}_{c}(\hat{u}_{1:N}, \lambda_{1:N+1}^{*}) \\
        & + \delta \underbar{$u$}_{1:N}^{\top} \nabla \mathcal{L}_{c}(\hat{u}_{1:N}, \lambda_{1:N+1}^{*}) \\
        & + \frac{1}{2} \, \delta \underbar{$u$}_{1:N}^{\top} \nabla^2 \mathcal{L}_{c}(\hat{u}_{1:N}, \lambda_{1:N+1}^{*}) \, \delta \underbar{$u$}_{1:N},
    \end{aligned}
\end{equation}
where $\delta \underbar{$u$}_{1:N}^{\top} = \left[ \delta u_1^\top, \hdots , \delta u_N^\top \right]^\top$ and $\delta u_t = u_t - \hat{u}_t, \forall \, 1 \leq t \leq N$. Additionally, we define the state deviations $\delta \underbar{$x$}_{1:N+1} = \left[ \delta x_1^\top, \hdots, \delta x_{N+1}^\top\ \right]^\top$ where $\delta x_t = x_t - \hat{x}_t, \forall \, 1 \leq t \leq N+1$. The symbols $\nabla \mathcal{L}_{c}(\cdot)$ and $\nabla^2 \mathcal{L}_{c}(\cdot)$ stand for the gradient and Hessian of $\mathcal{L}_{c}(\cdot)$, respectively. Following \cite{dunn1989efficient}, we expand the linear and quadratic terms of $\mathcal{L}_{c}(\cdot)$ as follows:
\begin{align}
    \delta \underbar{$u$}_{1:N}^{\top} \nabla \mathcal{L}_{c} & = \sum_{t=1}^{N} \delta u_{t}^{\top} d_{t}, \\[-0.5em]
    \delta \underbar{$u$}_{1:N}^{\top} \nabla^2 \mathcal{L}_{c} \, \delta \underbar{$u$}_{1:N} & = 
    \begin{aligned}[t]
        & \delta x_{N+1}^\top P_{N+1} \, \delta x_{N+1} + \sum_{t=1}^{N}\delta x_t^\top P_t \, \delta x_t \\[-0.5em]
        & + \sum_{t=1}^{N}\delta u_t^\top R_t \, \delta u_t + 2 \sum_{t=1}^{N} \delta x_t^\top M_t \, \delta u_t,\\[-1em]
    \end{aligned}
\end{align}
where we have defined $P_{N+1} = \partial^2 l_{N+1} / \partial x_{N+1}\partial x_{N+1}$ and 
\begin{align}
    P_t & = \frac{\partial^2 H_t}{\partial x_t\partial x_t} = \frac{\partial^2 l_t}{\partial x_t \partial x_t} + \frac{\partial^2 c_t}{\partial x_t \partial x_t} + \Big[ \lambda_{t+1}^{*} \Big]^\top \cdot \frac{\partial^2 f_t}{\partial x_t \partial x_t}, \\
    R_t & = \frac{\partial^2 H_t}{\partial u_t\partial u_t} = \frac{\partial^2 l_t}{\partial u_t \partial u_t} + \frac{\partial^2 c_t}{\partial u_t \partial u_t} + \Big[ \lambda_{t+1}^{*} \Big]^\top \cdot \frac{\partial^2 f_t}{\partial u_t \partial u_t},\\
    M_t & = \frac{\partial^2 H_t}{\partial x_t\partial u_t} = \frac{\partial^2 l_t}{\partial x_t \partial u_t} +  \frac{\partial^2 c_t}{\partial x_t \partial u_t} + \Big[ \lambda_{t+1}^{*} \Big]^\top \cdot \frac{\partial^2 f_t}{\partial x_t \partial u_t},\\
    d_t & = \frac{\partial H_t}{\partial u_t} = \frac{\partial l_t}{\partial u_t} + \frac{\partial c_t}{\partial u_t} + \Big[ \lambda_{t+1}^{*} \Big]^\top \frac{\partial f_t}{\partial u_t}.
\end{align}
Here, the notation $(A \cdot b)$ refers to a tensor contraction operation so that $(A \cdot b)_{ij} = \sum_{k} A_{ijk} \, b_{k}$. Additionally, all derivatives are evaluated at the nominal trajectory $\{ \hat{x}_{1:N+1}, \hat{u}_{1:N} \}$.

Given the expansion of $\mathcal{L}_{c}(\cdot)$ in~\eqref{eq:quad-lagrangian}, a Newton iteration on $\delta \underbar{$u$}_{1:N}$ is characterized by the minimization problem
\begin{equation}
    \label{eq:min-quad-lagrangian}
    \min_{\delta \underbar{$u$}_{1:N}}
    \begin{aligned}[t]
        & \frac{1}{2} \, \delta \underbar{$u$}_{1:N}^{\top} \nabla^2 \mathcal{L}_{c}(\hat{u}_{1:N}, \lambda_{1:N+1}^{*}) \, \delta \underbar{$u$}_{1:N} \\
        &+ \delta \underbar{$u$}_{1:N}^{\top} \nabla \mathcal{L}_{c}(\hat{u}_{1:N}, \lambda_{1:N+1}^{*}) + \alpha \, || \delta \underbar{$u$}_{1:N} ||_2^2.
    \end{aligned}
\end{equation}
Note that we have augmented the objective with a regularizing term in the style of a Levenberg--Marquardt algorithm with $\alpha \geq 0$ to guarantee a valid descent direction \cite{nocedal1999numerical}. After substituting the gradient and Hessian terms, we find the objective \eqref{eq:min-quad-lagrangian} to resemble a linear quadratic control problem over an augmented cost, subject to differential dynamics equations governing $\delta x_{t+1}$ as a function of $\delta x_{t}$ and $\delta u_{t}$ \cite{dunn1989efficient}
\begin{equation}
    \label{eq:linearized-ocp}
    \begin{split}
        &\min_{\delta \underbar{$u$}_{1:N}} \psi_{N+1}(\delta x_{N+1}) + \sum_{t=1}^{N} \psi_{t}(\delta x_{t}, \delta u_{t})\\
        &\textrm{subject to}\\
        &\delta x_{t+1} = \frac{\partial f_t}{\partial x_t} \, \delta x_t + \frac{\partial f_t}{\partial u_t} \, \delta u_t, \,\, \delta x_{1} = 0,
    \end{split}
\end{equation}
where the augmented cost terms are given by
\begin{align}
    \label{eq:linearized-ocp-param}
    \psi_{N+1}(\delta x_{N+1}) & = \frac{1}{2} \, \delta x_{N+1}^\top P_{N+1} \, \delta x_{N+1}, \\
    \psi_{t}(\delta x_{t}, \delta u_{t}) & = 
    \begin{aligned}[t]
        & \frac{1}{2} \delta x_t^\top P_t \, \delta x_t + \delta x_t^\top M_t \, \delta u_t \\
        & + \frac{1}{2} \delta u_t^\top (R_t + \alpha  I) \, \delta u_t + \delta u_t^\top d_t, 
    \end{aligned}
\end{align}
and $I\in\mathbb{R}^{d_u \times d_u}$ is the identity matrix. The solution to this differential control problem is obtained by solving a discrete-time Bellman equation \cite{lewis2012optimal}. This procedure requires defining a state-value function at time step $t$ as the cost-to-go
\begin{equation}
    \label{eq:linearized-ocp-bellman}
    V_t^{*}(\delta x_t) = \min_{\delta u_{t}} \psi_t(\delta x_t, \delta u_t) + V_{t+1}^{*}(\delta x_{t+1}),
\end{equation}
where $V_{t}^{*}(\delta x_{t})$, $1 \leq t \leq N+1$, maintains a quadratic form
\begin{equation}
    V_t^{*}(\delta x_t) = \frac{1}{2} \, \delta x_t^\top S_t \, \delta x_t + \delta x_t^\top s_t + \mathrm{const},
\end{equation}
starting from the boundary condition at $t=N+1$
\begin{equation}
    V_{N+1}^{*}(\delta x_{N+1}) = \frac{1}{2} \, \delta x_{N+1}^\top P_{N+1} \, \delta x_{N+1}.
\end{equation}
Computing the value function parameters $S_t$ and $s_t$ for all time steps via a sequential Riccati recursion incurs linear complexity in the length of the horizon $\mathcal{O}(N)$. It is possible to perform the same computation in logarithmic time $\mathcal{O}(\log N)$ by relying on parallel prefix-sums as in \cite{sarkka2022temporal}. In the following, we apply an analogous approach to our modified objective~\eqref{eq:linearized-ocp} and its corresponding Bellman equation~\eqref{eq:linearized-ocp-bellman}.
\begin{algorithm}[t]
    \caption{Parallel Value Function Pass}
    \label{alg:value-fun-parallel}
    \begin{algorithmic}[1]
        \State \textbf{Input}: Co-states: $\lambda_{1:N+1}^{*}$, nominal: $\{ \hat{x}_{1:N+1}, \hat{u}_{1:N} \}$.
        \State Evaluate partial derivatives of $l_{1:N+1}(\cdot)$, $c_{1:N}(\cdot)$, and $f_{1:N}(\cdot)$ at $\{ \hat{x}_{1:N+1}, \hat{u}_{1:N} \}$ required for~\eqref{eq:linearized-ocp} and \eqref{eq:linearized-ocp-param}.
        \State Initialize associative elements $\nu_{1:N+1}$ according to \eqref{eq:value-fun-elem-init}.
        \State Define the associative operator $\otimes$ as in \eqref{eq:value-fun-operator}
        \State Execute a parallel associative scan to compute $V_{1, t}(\delta x_1, \delta x_t) = \nu_1 \otimes \nu_2 \otimes \cdots \otimes \nu_{t-1}$ and corresponding $V_t^{*}(\delta x_t) = \nu_{t} \otimes \cdots \otimes \nu_{N+1}$, for $t= 1, \ldots, N+1$.
        \State Obtain value function parameters $S_t \gets Y_{t,\triangleleft}, s_t \gets -\eta_{t,\triangleleft}$, for $t= 1, \ldots, N+1$.
        \State Compute control laws $\Gamma_{1:N}, \gamma_{1:N}$ according to \eqref{eq:value-fun-feedback}.
        \State \textbf{Output}: $\Gamma_{1:N}, \gamma_{1:N}$.
    \end{algorithmic}
\end{algorithm}
Following the recipe in \cite{sarkka2022temporal}, we can construct an associative operator over our differential value functions from~\eqref{eq:linearized-ocp-bellman}
\begin{align}
    \label{eq:value-fun-operator}
    V_{j,i}(\delta x_j, \delta x_i) & = V_{j,t}(\delta x_j, \delta x_t) \otimes V_{t,i}(\delta x_t, \delta x_i) \\
    & = \min_{\delta x_t} \left[ V_{j,t}(\delta x_j, \delta x_t) + V_{t,i}(\delta x_t, \delta x_i) \right],
\end{align}
where an element $V_{j,i}(\delta x_j, \delta x_i)$ is defined as the \emph{conditional} value function that accumulates cost between two time steps $j \rightarrow i$ where $i > j$ and
\begin{equation}
    \label{eq:value-fun-def}
    \begin{split}
        V_{j, i}(\delta x_j, \delta x_i)
    = &\min_{\delta u_{j:i-1}} \sum_{t=j}^{i-1} \psi_{t}(\delta x_t, \delta u_t) \\
    &\textrm{subject to}\\
     &\delta x_{t+1} = \frac{\partial f_t}{\partial x_t} \, \delta x_t + \frac{\partial f_t}{\partial u_t} \, \delta u_t.
    \end{split}
\end{equation}
Consequently, an element $V_{j,i}(\delta x_j, \delta x_i)$ can be defined in terms of the convex -- and tractable -- dual problem of~\eqref{eq:value-fun-def}
\begin{equation}
    \label{eq:value-fun-elements}
    \begin{split}
        V_{j,i}(\delta x_j, \delta x_i) = 
        \max_{p}
        \begin{aligned}[t]
            & \frac{1}{2} \, \delta x_j^\top Y_{j,i} \, \delta x_{j} \\
            & - \delta x_j \, \eta_{j,i} - \frac{1}{2} p^{\top} C_{j,i} \, p \\
            & - p^{\top} \, \Big[ \delta x_i - A_{j,i} \, \delta x_j - b_{j,i} \Big],
        \end{aligned}
    \end{split}
\end{equation}
where $p \in \mathbb{R}^{d_x}$ is the Lagrangian dual variable. The variables $Y_{j,i}, C_{j,i}, A_{j,i} \in \mathbb{R}^{d_x \times d_x}$ and $\eta_{j,i}, b_{j,i} \in \mathbb{R}^{d_x}$, which parameterize the dual form of an element $V_{j,i}$ in~\eqref{eq:value-fun-elements}, arise from applying the  associated operator to previous elements $V_{j,t}$ and $V_{t,i}$ as in \eqref{eq:value-fun-operator}, and are given by
\begin{equation}
    \label{eq:value-fun-param}
    \begin{aligned}
        A_{j,i} & = A_{t,i} \, (I + C_{j,t} \, Y_{t,i})^{-1} A_{j,t}, \\
        Y_{j,i} & = A_{j,t}^\top (I + Y_{t,i} \, C_{j,t})^{-1} Y_{t,i} \, A_{j,t} + Y_{j,t}, \\
        C_{j,i} & = A_{t, i} \, (I + C_{j, t} \, Y_{t,i})^{-1} C_{j,t} \, A_{t, i}^\top + C_{t,i}, \\
        \eta_{j,i} & = A_{j,t}^\top \, (I + Y_{t,i} \, C_{j,t})^{-1} (\eta_{t,i} - Y_{t,i} \, b_{j,t}) + \eta_{j,t}, \\
        b_{j,i} & = A_{t,i} \, (I + C_{j,t} \, Y_{t,i})^{-1} (b_{j,t} + C_{j,t} \, \eta_{t,i}) + b_{j,i}.
    \end{aligned}
\end{equation}
We do not provide a full derivation of these formulas due to space constraints. Instead, we refer the reader to the general scheme from \cite{sarkka2022temporal}, as our approach follows similar steps. 

Analogous to Section~\ref{sec:co-states-scan}, the associative  parameters~\eqref{eq:value-fun-param} are initialized at consecutive time steps $t\rightarrow t+1$ as follows
\begin{equation}
    \label{eq:value-fun-elem-init}
    \nu_t = V_{t, t+1}(\delta x_t, \delta x_{t+1}),  \quad \forall \, 1 \leq t \leq N+1,
\end{equation}
and their parameters corresponding to \eqref{eq:value-fun-param} are initialized to
\begin{align}
    A_{t, t+1} & = \frac{\partial f_t}{\partial x_t} - \frac{\partial f_t}{\partial u_t} \, \tilde{R}_t^{-1} \, M_t^\top, \\
    Y_{t, t+1} & = P_t - M_t \, \tilde{R}_t^{-1} \, M_t^\top, \\
    C_{t, t+1} & = \frac{\partial f_t}{\partial u_t} \, \tilde{R}_t^{-1} \left[\frac{\partial f_t}{\partial u_t}\right]^\top, \\
    \eta_{t, t+1} & = \Big[ P_t - M_t \, \tilde{R}_t^{-1} \, M_t^\top \Big] r_t,\\
    b_{t, t+1} & = \frac{\partial f_t}{\partial u_t} \, \tilde{R}_t^{-1} \, M_t^\top r_t + \frac{\partial f_t}{\partial u_t} \, q_t,
\end{align}
with 
\begin{align}
    q_t & = - \Big[ \tilde{R}_t - M_t^\top P_t^{-1} M_t \Big]^{-1} d_t, \\
    r_t & = - P_t^{-1} M_t \, q_t, \quad \tilde{R}_{t} = R_t + \alpha I,
\end{align}
$\forall \, 1 \leq t \leq N$, whereas the parameters of the element $\nu_{N+1}$ are all initialized to zero except for $Y_{N+1,\triangleleft} = P_{N+1}$, where the symbol $\triangleleft$ stands for the boundary index. 

Now that we have defined the associative operator over value functions~\eqref{eq:value-fun-operator} and its elements~\eqref{eq:value-fun-elements}, the cost of computing the recursion~\eqref{eq:linearized-ocp-bellman} is reduced to logarithmic time $\mathcal{O}(\log N)$ by performing the associative scans in parallel to compute
\begin{align}
    \label{eq:value-fun-prefix-sum}
    \nu_1 \otimes \nu_2 \otimes \cdots \otimes \nu_{t-1} & = V_{1,t}(\delta x_1, \delta x_t), \\
    \nu_{t} \otimes \cdots \otimes \nu_{N+1} & = V_t^{*}(\delta x_t).
\end{align}
The value function parameters are then obtained as
\begin{equation}
    S_t = Y_{t,\triangleleft}, \quad s_t = - \eta_{t,\triangleleft},
\end{equation}
based on which we can derive an optimal control deviation vector and an affine control law $\delta u_{t}^{*} = \Gamma_t \, \delta x_{t} + \gamma_t$, where
\begin{equation}
    \label{eq:value-fun-feedback}
    \begin{split}
    \gamma_t & = - Q_{t}^{-1} d_t - Q_{t}^{-1} \left[ \frac{\partial f_t}{\partial u_t} \right]^\top \! s_{t+1}, \\
    \Gamma_t & = - Q_{t}^{-1} M_t -Q_{t}^{-1} \left[\frac{\partial f_t}{\partial u_t}\right]^\top \! S_{t+1} \, \frac{\partial f_t}{\partial x_t}, \\
    Q_{t} & = \tilde{R}_t + \left[\frac{\partial f_t}{\partial u_t}\right]^\top \! S_{t+1} \, \frac{\partial f_t}{\partial u_t}.
    \end{split}
\end{equation}
Algorithm~\ref{alg:value-fun-parallel} summarizes the parallel value function pass.

\subsection{Parallel State Propagation Pass}
The previous sections derived parallel-in-time algorithms for computing the optimal co-states $\lambda_{N+1}^{*}$, value functions $V_{1:N+1}^{*}$, and control deviation vector $\delta u_{1:N}^{*}$ conditioned on a nominal trajectory $\{ \hat{x}_{1:N+1}, \hat{x}_{1:N} \}$. These computations constitute a single iteration of Newton's algorithm. In order to perform the full iterative procedure, we need a routine to update the nominal trajectory. In this section, we derive a parallel associative scan, to compute the optimal state deviation vectors $\delta x_{1:N+1}^{*}$ that correspond to $\delta u_{1:N}^{*}$, according to \eqref{eq:linearized-ocp}, in logarithmic time $\mathcal{O}(\log N)$. The approach is analogous to scan from Section~\ref{sec:co-states-scan}, albeit defined in a reverse direction.

We start by defining the closed-loop differential dynamics
\begin{equation}
    \label{eq:closed-loop-dynamics}
    \delta x_{t+1}^{*} = \varphi_t (\delta x_t) = F_t \, \delta x_t^{*} + e_t,
\end{equation}
where $\delta x_1^{*} = 0$ and we assume
\begin{equation}
  F_t = \frac{\partial f_t}{\partial x_t} + \frac{\partial f_t}{\partial u_t} \, \Gamma_t, \quad e_t = \frac{\partial f_t}{\partial u_t} \, \gamma_t.
\end{equation}
The recursion in~\eqref{eq:closed-loop-dynamics} can be rewritten in terms of $\delta x_1^{*}$ using the function composition:
\begin{equation}
    \delta x_t^{*} = (\varphi_{t-1} \circ \cdots \circ \varphi_{1}) (\delta x_{1}^{*}).
\end{equation}
Next, let us define the function composition $\circ$ as the associative operator $\otimes$ over the state deviation elements
\begin{equation}
    \label{eq:state-prop-operator}
    \varphi_{j,i} = \varphi_{t,i} \circ \varphi_{j,t},
\end{equation}
where an element $\varphi_{j,i}$ results from applying the associative operator 
between the time steps $j \rightarrow i$ where $i > j$
\begin{equation}
    \label{eq:state-prop-elements}
    \begin{split}
        \varphi_{j,i}
        & = \left(\varphi_{i-1} \circ \cdots  \circ \varphi_j \right) (\delta x_j^{*}) \\
        & = F_{j,i} \, \delta x_j^{*} + e_{j,i}.
    \end{split}
\end{equation}
The parameters $F_{j,i}$ and $e_{j,i}$ are obtained after applying the associative operator defined in~\eqref{eq:state-prop-operator}  on the elements  $\varphi_{t,i}$ and $\varphi_{j,t}$ as follows
\begin{equation}
    \label{eq:state-prop-comb-rule}
    F_{j,i} = F_{t,i} \, F_{j,t}, \quad e_{j,i} = F_{t,i} \, e_{j,t} + e_{t,i}.
\end{equation}
We initialize the elements of the associative scan as $\vartheta_t = \varphi_{t,t+1}$ which can be written in terms of parameters as
\begin{equation}
    \label{eq:state-prop-elem-init}
    F_{t,t+1} = F_t, \quad e_{t,t+1} = e_t,
\end{equation}
for $1 < t \leq N$, while the initial element $\vartheta_1$ has parameters
\begin{equation}
    \label{eq:state-prop-init-elem-init}
    F_{1,2} = 0, \quad e_{1,2} = F_1 \, \delta x_1^{*} + e_1.
\end{equation}
All derivatives are evaluated with respect to the nominal trajectory $\{ \hat{x}_{1:N+1}, \hat{u}_{1:N} \}$. Using the associative operator~\eqref{eq:state-prop-operator} and its elements~\eqref{eq:state-prop-elements}, the propagation of the differential state equation to compute $\delta x_{1:N+1}^{*}$ is performed in logarithmic time $\mathcal{O}(\log N)$ via a parallel associative scan for computing
\begin{equation}
    \delta x_t^{*} = \vartheta_1 \otimes \vartheta_2 \otimes \hdots \otimes \vartheta_{t-1}, \quad \forall \, 1 \leq t \leq N.
\end{equation}
Algorithm~\ref{alg:state-prop-parallel} summarizes the parallel state propagation pass.

\begin{algorithm}[t]
    \caption{Parallel State Propagation Pass}
    \label{alg:state-prop-parallel}
    \begin{algorithmic}[1]
        \State \textbf{Input}: Controller: $\Gamma_{1:N}, \gamma_{1:N}$, nominal: $\{ x_{1:N+1}^{(n)}, u_{1:N}^{(n)} \}$.
        \State Initialize the initial associative element $\vartheta_1$ with \eqref{eq:state-prop-init-elem-init}.
        \State Initialize the associative elements $\vartheta_{2:N}$ with \eqref{eq:state-prop-elem-init}.
        \State Define the associative operator $\otimes$ as in \eqref{eq:state-prop-operator}.
        \State Execute a parallel associative scan to compute $\delta x_t^{*} = \vartheta_1 \otimes \vartheta_2 \otimes \cdots \otimes \vartheta_{t-1}$, for $t=1, \ldots, N$.
        \State Compute associated control deviations $\delta u_{1:N}^{*}$.
        \State \textbf{Output}: $\delta x_{1:N+1}^{*}, \delta u_{1:N}^{*}$.
    \end{algorithmic}
\end{algorithm}

\begin{algorithm}[t]
    \caption{Parallel Newton’s Method
for Solving \eqref{eq:lb-min} or \eqref{eq:admm-min}}
    \label{alg:par-iterative-newton}
    \setstretch{1.25}
    \small
    \begin{algorithmic}[1]
        \item \textbf{Input:} Initial nominal: $\{ x_{1:N+1}^{(0)}, u_{1:N}^{(0)} \}$.
        \State $n \gets 0$.
        \While{not converged}
            \State \hspace{-12pt} $\lambda_{1:N+1}^{(n)} \! \gets \!\textbf{Algorithm~\ref{alg:co-state-parallel}} \big( \{ x_{1:N+1}^{(n)}, u_{1:N}^{(n)} \} \big)$.
            \State \hspace{-12pt} $\Gamma_{1:N}^{(n)}, \gamma_{1:N}^{(n)} \! \gets \! \textbf{Algorithm~\ref{alg:value-fun-parallel}} \big(\lambda_{1:N+1}^{(n)}, \{ x_{1:N+1}^{(n)}, u_{1:N}^{(n)} \} \big)$.
            \State \hspace{-12pt} $\delta x_{1:N+1}^{(n)}, \delta u_{1:N}^{(n)} \! \gets \! \textbf{Algorithm~\ref{alg:state-prop-parallel}} \big( \Gamma_{1:N}^{(n)}, \gamma_{1:N}^{(n)}, \{ x_{1:N+1}^{(n)}, u_{1:N}^{(n)} \} \big)$.
            \State \hspace{-12pt} $x_{1:N+1}^{(n+1)} \gets x_{1:N+1}^{(n)} + \delta x_{1:N+1}^{(n)}$.
            \State \hspace{-12pt} $u_{1:N}^{(n+1)} \gets u_{1:N}^{(n)} + \delta u_{1:N}^{(n)}$.
            \State \hspace{-12pt} $n \gets n + 1$.
        \EndWhile
        \State \textbf{Output}: $x_{1:N+1}^{*}, u_{1:N}^{*}$.
        \end{algorithmic}
\end{algorithm}

\subsection{Parallel Newton's Method for Optimal Control}
In this section, we unify the previously presented parallel scan algorithms into a single iterative parallel-in-time Newton’s method, outlined in Algorithm~\ref{alg:par-iterative-newton}. This method can be embedded within IP methods or ADMM, to solve problems such as \eqref{eq:lb-min} and \eqref{eq:admm-min}, within constrained nonlinear optimal control problems. The general framework of IP methods and ADMM is provided in Section \ref{sec:background}. For thorough implementation guidelines, see \cite{boyd2011distributed, boyd2004convex, nocedal1999numerical}.%

\subsection{Convergence and Feasibility}
In this section, we discuss the convergence of the introduced algorithms.

In the log-barrier method, solving barrier problems and decreasing the barrier parameter repeatedly, leads to a feasible and optimal solution \cite{boyd2004convex}. In ADMM, both the primal and dual residual are shown to converge \cite{boyd2011distributed}, leading to the feasibility and optimality of the solution.

The temporal parallelization of Newton's method for optimal control inherits the convergence properties of the sequential implementation \cite{dunn1989efficient, de1988differential}. The trajectory iterates resulting from the Newton updates converge to a local minimizer for any nearby initial candidate \cite{dunn1989efficient}. The regularizing parameter $\alpha$ introduced in \eqref{eq:min-quad-lagrangian} gets updated based on the trust region rule described in (2.21) of \cite{madsen2004methods} and ensures a decent direction similar to the Levenberg--Marquardt strategy \cite{nocedal1999numerical}.

The log-barrier method ensures recursive feasibility as long as the initial solution guess is feasible. The logarithmic barrier function enforces feasibility by penalizing any deviations from valid iterates \cite{boyd2004convex}. In contrast, ADMM does not face recursive feasibility concerns, as its optimization process can begin from an infeasible trajectory \cite{ma2022alternating}. The method's computational result enables longer-horizon MPC optimization, arguably enhancing stability.

\section{Numerical Results}
In this section, we implement the parallel associative scans, defined by their operators, elements, and initialization routines, using the JAX software package~\cite{jax2018github}. Numerical evaluations are performed on an NVIDIA Tesla A100 GPU. However, these parallel scan algorithms are not restricted to GPU architectures and can be executed on various types of massively parallel computing systems. 

Our evaluation emphasizes the advantages of our parallel-in-time version of Newton's method for nonlinear optimal control, in comparison to the sequential-in-time approach initially introduced in \cite{dunn1989efficient} as well as interior point differential dynamic programming (IPDDP) \cite{pavlov2021interior}. We assess the computational efficiency of our algorithm within iterated interior point and ADMM methods, as discussed in Section~\ref{sec:constrained-optimal-control}, both of which involve repeatedly solving nonlinear optimal control problems over numerous iterations, significantly compounding the computational burden. 

We evaluate our method on two control-constrained dynamical systems: a torque-constrained pendulum and a force-constrained cart-pole. In each case, we perform constrained trajectory optimization over a fixed time, with varying sampling frequencies, resulting in trajectory optimization problems with variable number of planning steps. Furthermore, we run MPC simulations as well.

The initial trajectory is obtained from a nonlinear dynamics rollout for a random, zero-mean normal distributed control sequence. The optimal trajectory resulting from an IP or ADMM iteration gets recycled as the initial trajectory for the next iteration. Further practical details of our implementation can be found in the linked code base\footnote{\url{https://github.com/casiacob/constrained-parallel-control}}.

\subsection{Torque-Limited Pendulum Swing-Up Problem}
In our first experiment, we tackle a pendulum swing-up and stabilization task with torque constraints \cite{tedrake2023underactuated}
\begin{equation}
        \ddot{\theta} = -\frac{g}{l}\,\sin\theta + \frac{\tau - b\, \dot{\theta}}{m\,l^2},
\end{equation}
where $\theta$ is the pendulum angle, $\tau$ is the input torque, $g=9.81\mathrm{m/s^2}$ is the gravitational acceleration, $l=1\mathrm{m}$ is the pendulum length, $m=1\mathrm{kg}$ is the pendulum mass, and $b=10^{-3} \mathrm{kg}\,{\mathrm{m}}/{\mathrm{s}}$ is the damping. The pendulum dynamics are integrated using Euler's method across different sampling frequencies, resulting in planning horizons ranging from $N = 20$ to $N = 1000$. The barrier parameter in the interior point method starts at $\mu = 0.1$ and is reduced by a factor of $\zeta = 0.2$ until $\mu \leq 10^{-4}$. For ADMM, we set the penalty parameter to $\rho = 1$, running the algorithm until the primal and dual residuals satisfy $||r_p||_\infty, ||r_d||_\infty \leq 10^{-2}$, where $||\cdot||_\infty$ represents the infinity norm of a vector. Both methods share a common inner convergence threshold in Algorithm~\ref{alg:par-iterative-newton}. Figure~\ref{fig:pendulum-runtime} shows the average runtime over $10$ runs for both the IP method and ADMM, comparing the performance of our parallel-in-time algorithm against the sequential-in-time method from~\cite{dunn1989efficient}, The results demonstrate the scalability advantage of our parallel-in-time approach as the planning horizon increases.

\subsection{Force-Limited Cart-Pole Swing-Up Problem}
In the second experiment, we simulate a force-limited cart-pole swing-up and stabilization task \cite{tedrake2023underactuated}
\begin{align}
    \ddot{p} &= \frac{1}{m_c + m_p\,\sin^2\theta}(F + m_p\,\sin\theta\,(l\,\dot{\theta} + g\,\cos\theta)),\\
    \ddot{\theta} &= 
    \begin{aligned}[t]
    \frac{1}{l(m_c+m_p\,\sin^2\theta)}(&-F\,\cos\theta - m_p\,l\,\dot{\theta}^2\cos\theta\,\sin\theta\\
    &- (m_c+m_p)\,g\sin\theta),
    \end{aligned}
\end{align}
where $p$ is the cart position, $\theta$ is the pole angle, $F$ is the input force pulling the cart, $g=9.81 \mathrm{m}/\mathrm{s^2}$ is the gravitational acceleration, $l=0.5 \mathrm{m}$ is the pole length, $m_c=10\mathrm{kg}$ is the cart mass and $m_p=1 \mathrm{kg}$ is the pole mass.
Similar to the pendulum experiment, we integrate the cart-pole dynamics using Euler's method at varying frequencies, producing planning horizons in the range $N = 20$ to $N = 1000$. For the interior point method, the barrier parameter is initialized at $\mu = 0.1$ and decreased by $\zeta = 0.2$ until it reaches $\mu \leq 10^{-4}$. In the ADMM method, the penalty parameter is set to $\rho = 0.5$, and the algorithm continues until both the primal and dual residuals meet the criteria $||r_p||_\infty ,||r_d||_\infty \leq 10^{-2}$. As before, Algorithm~\ref{alg:par-iterative-newton} applies the same convergence threshold for both methods. Figure~\ref{fig:cartpole-runtime} illustrates the mean runtime across $10$ runs, comparing our parallel-in-time algorithm with the sequential-in-time method from~\cite{dunn1989efficient}. Similar to the first experiment, the parallel-in-time method proves more computationally efficient, scaling more gracefully for longer planning horizons.
\subsection{Model Predictive Control Simulation}
Finally, we run a model predictive control simulation for both systems. We target vertical angle stabilization for the pendulum and the cart-pole and position control for the cart. The control loops are simulated for 4 seconds. The dynamics are sampled using Euler's method, at a rate of 100Hz. The optimization horizons are fixed to $N=60$ steps. Our method stabilizes both systems while satisfying the control constraints, $|\tau| < 5\mathrm{Nm}$, and $|F|<60 \mathrm{N}$. The resulting trajectories are displayed in Figure
\ref{fig:cartpole-mpc}.

\section{Conclusion}
We presented a highly efficient parallel realization of Newton's method for dynamic optimization in constrained optimal control problems. Our approach introduces an algorithmic framework that harnesses the power of massively parallel computation to reduce the computational burden associated with long planning horizons significantly. By deploying our method on GPUs and evaluating it across a series of simulated benchmarks, we demonstrated its advantages in performance and scalability compared to traditional recursive approaches. This highlights the potential of our method to accelerate large-scale optimization tasks in nonlinear control systems.
\begin{figure}[t]
    \centering
\begin{tikzpicture}

\definecolor{darkgray176}{RGB}{176,176,176}
\definecolor{darkorange25512714}{RGB}{255,127,14}
\definecolor{lightgray204}{RGB}{204,204,204}
\definecolor{steelblue31119180}{RGB}{31,119,180}

    \begin{groupplot}[
        group style={
            group name=my plots,
            group size=2 by 1, 
            horizontal sep=5pt
        }, 
        width=5.0cm, 
        height=4.5cm, 
        grid=both,
        minor x tick num=3,
        grid style={line width=.1pt, draw=gray!20},
        major grid style={line width=.1pt, draw=gray!60},
        ytick distance=10^1,
        ytick={0.1,1,10,100,1000},
        yticklabels={$10^{-1}$,,$10^{1}$,,$10^{3}$},
        ymin=0.07, ymax=1100,
        ylabel shift=-0.25cm,
    ]
    \nextgroupplot[
        title={Interior Point},
        xmode=log,
        ymode=log,
        xlabel={Horizon $N$},
        ylabel={Runtime [s]},
        xtick align=inside,
        ytick align=inside,
        legend entries={Parallel;,Sequential;,IPDDP},
        legend to name=common_legend,
        legend style={legend columns=3},
    ]
        \addplot [
            semithick, 
            black, 
            mark=o, 
            mark size=2, 
            mark options={solid}
        ] table {%
            20 0.165869498252869
            40 0.210235667228699
            80 0.254219579696655
            100 0.258720445632935
            200 0.290592741966248
            400 0.340044283866882
            800 0.397438740730286
            1000 0.399478697776794

        };
        \addplot [
            semithick, 
            black, 
            mark=triangle, 
            mark size=2, 
            mark options={solid}
        ] table {%
            20 0.259859371185303
            40 0.505869269371033
            80  1.05387041568756
            100 1.33567306995392
            200 2.64988203048706
            400 5.33660573959351
            800 10.6862892627716
            1000 13.2525683403015
        };

        \addplot [
            semithick, 
            black, 
            mark=square, 
            mark size=2, 
            mark options={solid}
        ] table {%
            20 0.15359962 
            40 0.26941674 
            80 0.5138206  
            100 0.69460416 
            200 1.39813712 
            400 2.72683363
            800 5.54038296 
            1000 6.77793486
        };

    \nextgroupplot[
        title={ADMM},
        xmode=log,
        ymode=log,
        xlabel={Horizon $N$},
        ylabel={},
        yticklabel={\empty},
        xtick align=inside,
        ytick align=inside,
        ymin=0.07, ymax=1100,
    ]
        \addplot [
            semithick, 
            black, 
            mark=o, 
            mark size=2, 
            mark options={solid}
        ] table {%
            20 0.599585127830505
            40  2.52651727199554
            80  3.98368370532989
            100  4.3857483625412
            200  6.25698399543762
            400  7.67760238647461
            800  8.2451893568039
            1000 8.44810690879822

        };


        \addplot [
            semithick, 
            black, 
            mark=triangle, 
            mark size=2, 
            mark options={solid}
        ] table {%
            20 1.01457891464233
            40 7.11924307346344
            80 18.5311008453369
            100 25.3128653287888
            200 60.4472321748734
            400 133.003935575485
            800 252.717141699791
            1000 324.930444955826
        };
    
\end{groupplot}
\node[anchor=north] (title-x) at ($(my plots c1r1.south east)!0.5!(my plots c2r1.south west)-(0,1.0cm)$) {\ref{common_legend}};
\end{tikzpicture}
    \vskip -10pt
    \caption{Runtime average of parallel and sequential realization of Newton-type ADMM and interior point methods applied to a torque-constrained pendulum swing-up problem solved for various planning horizons.}
    \label{fig:pendulum-runtime}
    \vspace{-0.4cm}
\end{figure}

\begin{figure}[t]
    \centering
\begin{tikzpicture}

\definecolor{darkgray176}{RGB}{176,176,176}
\definecolor{darkorange25512714}{RGB}{255,127,14}
\definecolor{lightgray204}{RGB}{204,204,204}
\definecolor{steelblue31119180}{RGB}{31,119,180}

    \begin{groupplot}[
        group style={
            group name=my plots,
            group size=2 by 1, 
            horizontal sep=5pt
        }, 
        width=5.0cm, 
        height=4.5cm, 
        grid=both,
        minor x tick num=3,
        grid style={line width=.1pt, draw=gray!20},
        major grid style={line width=.1pt, draw=gray!60},
        ytick distance=10^1,
        ytick={0.1,1,10,100,1000},
        yticklabels={$10^{-1}$,,$10^{1}$,,$10^{3}$},
        ymin=0.07, ymax=2000,
        ylabel shift=-0.25cm,
    ]
    \nextgroupplot[
        title={Interior Point},
        xmode=log,
        ymode=log,
        xlabel={Horizon $N$},
        ylabel={Runtime [s]},
        xtick align=inside,
        ytick align=inside,
        legend entries={Parallel;,Sequential;,IPDDP},
        legend to name=common_legend1,
        legend style={legend columns=3},
    ]
        \addplot [
            semithick, 
            black, 
            mark=o, 
            mark size=2, 
            mark options={solid}
        ] table {%
            20 0.27523365020752
            40 0.39802610874176
            80 0.524001717567444
            100 0.564911413192749
            200 0.754844188690185
            400 0.817452573776245
            800 0.874869179725647
            1000 0.974847149848938

        };


        \addplot [
            semithick, 
            black, 
            mark=triangle, 
            mark size=2, 
            mark options={solid}
        ] table {%
            20 0.431954693794251
            40 0.894152617454529
            80 1.95449202060699
            100 2.50097727775574
            200 4.95996022224426
            400 9.95127701759338
            800 20.1917254924774
            1000 26.7140133142471
        };

        \addplot [
            semithick, 
            black, 
            mark=square, 
            mark size=2, 
            mark options={solid}
        ] table {%
            20 0.34495893 
            40 0.57070098 
            80 1.60221007 
            100 1.11908929 
            200 1.7820322  
            400 6.34651108
            800 7.94645014 
            1000 9.69574673
        };

    \nextgroupplot[
        title={ADMM},
        xmode=log,
        ymode=log,
        xlabel={Horizon $N$},
        ylabel={},
        yticklabel={\empty},
        xtick align=inside,
        ytick align=inside,
        ymin=0.07, ymax=2000,
    ]
        \addplot [
            semithick, 
            black, 
            mark=o, 
            mark size=2, 
            mark options={solid}
        ] table {%
            20 11.6474112510681
            40 14.1542631387711
            80 17.3078180551529
            100 17.7634773254395
            200 21.0172792434692
            400 23.8359166383743
            800 26.7202964067459
            1000 27.0246778011322

        };


        \addplot [
            semithick, 
            black, 
            mark=triangle, 
            mark size=2, 
            mark options={solid}
        ] table {%
            20 19.5968326330185
            40 40.2055332183838
            80 85.999801492691
            100 104.240215945244
            200 212.635953617096
            400 424.69156460762
            800 869.105820465088
            1000 1094.50579841137

        };    

\end{groupplot}
\node[anchor=north] (title-x) at ($(my plots c1r1.south east)!0.5!(my plots c2r1.south west)-(0,1.0cm)$) {\ref{common_legend1}};
\end{tikzpicture}
    \vskip -10pt
    \caption{Runtime average of parallel and sequential realization of Newton-type ADMM and interior point methods applied to a force-constrained cart-pole swing-up problem solved for various planning horizons.}
    \label{fig:cartpole-runtime}
    \vskip 5pt
\end{figure}

\begin{figure}[t]
    \centering
    \input{figures/cartpole_mpc}
    \vskip -10pt
    \caption{Model predictive control of constrained pendulum (left) and constrained cart-pole (right).}
    \label{fig:cartpole-mpc}
    \vskip 5pt
\end{figure}
\vspace{-0.15cm}
\bibliographystyle{IEEEtran}
\bibliography{IEEEabrv,references}
\end{document}